\documentclass[aop,seceqn,MSNbibl,citesort,dvips]{arximspdf}
\usepackage{graphicx}

\doi{10.1214/10-AOP581}
\volume{39}
\issue{4}
\pubyear{2011}
\firstpage{1286}
\lastpage{1304}

\makeatletter

\newtheorem{theorem}{Theorem}
\newproclaim{definition}[theorem]{Definition}
\newtheorem{lemma}[theorem]{Lemma}
\newtheorem{proposition}[theorem]{Proposition}
\newproclaim{remark}[theorem]{Remark}

\newcommand{\PP}{\mathbb{P}}

\makeatother

\begin{document}
\begin{frontmatter}

\title{On monochromatic arm exponents for 2D critical percolation}
\runtitle{Monochromatic exponents for percolation}

\begin{aug}
\author[A]{\fnms{Vincent} \snm{Beffara}\corref{}\ead[label=e1]{Vincent.Beffara@ens-lyon.fr}} and
\author[B]{\fnms{Pierre} \snm{Nolin}\ead[label=e2]{nolin@cims.nyu.edu}\thanksref{t1}}
\thankstext{t1}{Supported by NSF Grant OISE-07-30136.}
\runauthor{V. Beffara and P. Nolin}
\affiliation{ENS Lyon and the Courant Institute}
\address[A]{UMPA---ENS Lyon\\
46, all\'ee d'Italie\\
F-69364 Lyon cedex 07\\
France\\
\printead{e1}}
\address[B]{Courant Institute of \\
\quad Mathematical Sciences\\
New York University\\
251 Mercer Street\\
New York, New York 10012\\
USA\\
\printead{e2}}
\end{aug}

% HISTORY:
\received{\smonth{6} \syear{2009}}
\revised{\smonth{7} \syear{2010}}

% ABSTRACT
%
\begin{abstract}
We investigate the so-called \textit{monochromatic arm exponents} for
critical percolation in two dimensions. These exponents, describing
the probability of observing $j$ disjoint macroscopic paths, are
shown to exist and to form a different family from the (now well
understood) polychromatic exponents. More specifically, our main
result is that the monochromatic $j$-arm exponent is strictly
between the polychromatic $j$-arm and $(j+1)$-arm exponents.
\end{abstract}

% KEYWORDS
%

\end{frontmatter}

%s1 ###
\section{Introduction}

Percolation is one of the most-studied discrete models in statistical
physics. The usual setup is that of \textit{bond percolation} on the
square lattice~$\mathbb Z^2$, where each bond is open (resp., closed)
with probability $p\in(0,1)$ (resp., $1-p$), independently of the
others. This model exhibits a phase transition at a critical point
$p_c\in(0,1)$ (in this particular case, $p_c=1/2$): for $p<p_c$, almost
surely all connected components are finite, while for $p>p_c$ there
exists a unique infinite component with density $\theta(p)>0$. Site
percolation is defined in a similar fashion, the difference being that
the vertices are open or closed, instead of the edges; one can then see
it as a random coloring of the lattice, and use the terms \textit{black}
and \textit{white} in place of \textit{open} and \textit{closed}.

The behavior of percolation away from the critical point is well
understood; however it is only recently that precise results have been
obtained at and near criticality. For critical site percolation on the
regular triangular lattice, the proof of conformal invariance in the
scaling limit was obtained by Smirnov~\cite{smirnovperco}, and
$\mathrm{SLE}$ processes, as introduced by Schramm~\cite{schrammlerw}
and further studied by Lawler, Schramm and Werner
\cite{wernervalue,wernervalue2}, provide an explicit description of
the interfaces (in the scaling limit) in terms of $\mathrm{SLE}(6)$
(see, e.g.,~\cite{wernerparkcity}).

This description allows for the derivation of the so-called
\textit{polychromatic arm exponents}~\cite{smirnovexps}, which describe
the probability of observing connections across annuli of large modulus
by disjoint connected paths of specified colors (with at least one arm
of each color), and also the derivation of the one-arm exponent
\cite{lswonearm}. Combined with Kesten's scaling relations
\cite{kestenscaling}, these exponents then provide the existence and
the values of most of the other critical exponents, like, for example,
the exponent $\beta=5/36$ associated with the density of the infinite
cluster, as $p \downarrow p_c$,
\[
\theta(p) = (p-p_c)^{5/36 + o(1)}.
\]

On the other hand, very little is known concerning the
\textit{monochromatic} arm exponents (i.e., with all the
connections of the same color---see below for a formal definition)
with more than one arm. Here, the $\mathrm{SLE}$ approach does not seem
to work, and, correspondingly, there is no universally established
conjecture for the values of those exponents. One notable exception,
however, is the two-arm monochromatic exponent, for which an interpretation
in terms of $\mathrm{SLE}(6)$ is proposed at the end of
\cite{lswonearm}---but again no explicit value has been computed.
That particular exponent is actually of physical interest: known as the
\textit{backbone exponent}, it describes the ``skeleton'' of a percolation
cluster. Even the existence of these exponents is not clear, as there
does not seem to be any direct sub-additivity argument.

In this paper, we prove that the monochromatic exponents do exist, and
investigate how they are related to the polychromatic exponents. We show
that they have different values than their polychromatic counterparts.
As an illustration, our result implies that the backbone of a typical
large percolation cluster at criticality is much ``thinner'' than its
boundary.

%s2 ###
\section{Background}

%s2.1 ###
\subsection{The setting}

We restrict ourselves here to site percolation on the triangular
lattice, at criticality ($p=p_c=1/2$). Recall that it can be obtained by
coloring the faces of the honeycomb lattice randomly, each cell being
black or white with probability $1/2$ independently of the others. In
the following, we denote by $\PP= \PP_{1/2}$ the corresponding
probability measure on the set of configurations. Let us mention, however,
that many of the results of combinatorial nature based on
Russo--Seymour--Welsh-type estimates should also hold for bond percolation
on~$\mathbb{Z}^2$, due to the self-duality property of this lattice.

Let $S_n$ denote the ball of radius $n$ in the triangular lattice
(i.e., the intersection of the triangular lattice with the
Euclidean disc of radius $n$, though the specifics of the definition are
of little relevance), seen as a set of vertices. We will denote by
$\partial^i S_n$ (resp., $\partial^e S_n$) its internal (resp.,
external) boundary, that is, the set of vertices in (resp.,
outside)\vadjust{\goodbreak}
$S_n$ that have at least one neighbor outside (resp., in) $S_n$, and,
for $n<N$, by
\[
S_{n,N} := S_N \setminus S_n
\]
the annulus of radii $n$
and $N$. To describe critical and near-critical percolation, certain
exceptional events play a central role: the arm events, referring to the
existence of a number of crossings (``arms'') of $S_{n,N}$, the color of
each crossing (black or white) being prescribed.
\begin{definition}
Let $j \ge1$ be an integer and $\sigma= (\sigma_1,\ldots,\sigma_j)$
be a sequence of colors (black or white). For any two positive
integers $n < N$, a \textit{$(j,\sigma)$-arm configuration} in the
annulus $S_{n,N}$ is the data of $j$ disjoint monochromatic,
nonself-intersecting paths $(r_i)_{1 \leq i \leq j}$---the
\textit{arms}---connecting the inner boundary $\partial^e S_n$ and the
outer boundary $\partial^i S_N$ of the annulus, ordered
counterclockwise in a cyclic way, where the color of the arm $r_i$ is
given by $\sigma_i$. We denote by
%
%e2.1 ###
%
\begin{equation}
A_{j,\sigma}(n,N) := \bigl\{ \partial^e S_n \mathop{\leadsto}
_{j,\sigma} \partial^i S_N \bigr\}
\end{equation}
the corresponding event. It depends only on the state of the vertices
in $S_{n,N}$.
\end{definition}

We will write down color sequences by abbreviating colors, using $B$ and
$W$ for black and white, respectively. To avoid the obvious combinatorial
obstructions, we will also use the notation $n_0 = n_0(j)$ for the
smallest integer such that $j$ arms can possibly arrive on $\partial^e
S_{n_0}$ [$n_0(j)$ is of the order of $j$] and only consider annuli of
inner radius larger than $n_0$. This restriction will be done implicitly
in what follows.

The so-called \textit{color exchange trick} (noticed in
\cite{aharonyexponents,smirnovexps}) shows that for a fixed number
$j$ of arms, prescribing the color sequence $\sigma$ changes the
probability only by at most a constant factor, \textit{as long as both
colors are present in $\sigma$} (because an interface is needed to
proceed). The asymptotic behavior of that probability can be described
precisely using $\mathrm{SLE}(6)$: it is possible to prove the existence
of the (polychromatic) arm exponents and to derive their values
\cite{smirnovexps}, which had been predicted in the physics literature
(see, e.g.,~\cite{aharonyexponents} and the references therein).
\begin{theorem}
\label{thm:armcrit}
Fix $j \geq2$. Then for any color sequence $\sigma$ containing both
colors,
%
%e2.2 ###
%
\begin{equation}
\PP(A_{j,\sigma}(n,N)) = N^{-\alpha_j + o(1)}
\end{equation}
as $N \to\infty$ for any fixed $n \geq n_0(j)$, where $\alpha_j =
(j^2-1)/12$.
\end{theorem}

The value of the exponent for $j=1$ (corresponding to the probability of
observing one arm crossing the annulus) has also been established
\cite{lswonearm} and it is equal to $5/48$ (oddly enough formally
corresponding to $j=3/2$ in the above formula).\vadjust{\goodbreak}

For future reference, let us mention the following facts about critical
percolation that we will use.
\begin{enumerate}
\item\hypertarget{apriori}\mbox{} \textit{A priori bound for arm configurations}:
there exist constants $C, \varepsilon>0$ such that for all $n < N$,
%
%e2.3 ###
%
\begin{equation}
\PP(A_{1,B}(n,N)) = \PP(A_{1,W}(n,N)) \leq C
\biggl( \frac{n}{N} \biggr)^{\varepsilon}.
\end{equation}
For all $j\geq1$, there exist constants $c_j, \beta_j > 0$ such that
for all $n<N$,
%
%e2.4 ###
%
\begin{equation}
\label{apriorij}
\PP(A_{j,B\ldots BB}(n,N)) \geq c_j \biggl( \frac{n}{N}
\biggr)^{\beta_j}.
\end{equation}

\item\hypertarget{quasi_mult}\mbox{} \textit{Quasi-multiplicativity property}: For
any $j \geq1$ and any sequence $\sigma$, there exist constants
$C_1,C_2 > 0$ such that for all $n_1 < n_2 < n_3$,
\begin{eqnarray*}
&&
C_1 \PP(A_{j,\sigma}(n_1,n_2))
\PP(A_{j,\sigma} (n_2, n_3)) \\
&&\qquad\leq
\PP(A_{j,\sigma}(n_1,n_3))\\
&&\qquad\leq C_2 \PP(A_{j,\sigma}(n_1,n_2))
\PP(A_{j,\sigma}(n_2,n_3)).
\end{eqnarray*}
\end{enumerate}

The first of these two properties actually relies on the so-called
Russo--Seymour--Welsh (RSW) lower bounds, that we will use extensively in
various situations: roughly speaking, these bounds state that the
probability of crossing a given shape of fixed aspect ratio is bounded
below independently of the scale. For instance, the probability of
crossing a $3n \times n$ rectangle in its longer direction is bounded
below, uniformly as $n\to\infty$. We refer the reader
to~\cite{grimmettbook} for more details.

In the second one, the independence of arm events in disjoint annuli
implies that one can actually take $C_2=1$, which we will do from now
on. The lower bound is obtained using a so-called \textit{separation
lemma}, as first proved by Kesten~\cite{kestenscaling}, Lem\-ma~6
(technically, he does the proof in the case of four arms, but the
argument extends easily to the general case). The monochromatic case is
rather easier, as is follows directly from RSW estimates and Harris's
lemma.

%s2.2 ###
\subsection{A correlation inequality}

For two increasing events, the probability of their disjoint occurrence
can be bounded below by the classic van den Berg--Kesten (BK) inequality
\cite{vdBKBK}; Reimer's inequality, conjectured in~\cite{vdBKBK} and
proved in~\cite{reimerR}, extends it to the case of arbitrary events.
A key ingredient in our proof will be a not-that-classic correlation
inequality which is an intermediate step in the proof of Reimer's
inequality. Note that actually, in the simpler case of increasing events
(which is the only one we need here), this inequality was obtained
earlier by Talagrand~\cite{talagrandbk}; but we choose to follow the
more established name of the general result even for this particular
case.\vadjust{\goodbreak}

Instead of using the terminology in Reimer's original
paper~\cite{reimerR} we follow the rephrasing (with more ``probabilistic''
notation) of his proof in the review paper~\cite{BCRBKR}. Consider an
integer $n$, and $\Omega= \{0,1\}^n$. For any configuration $\omega
\in\Omega$ and any set of indices $S \subseteq\{1,\ldots,n\}$, we
introduce the cylinder
\[
[\omega]_S := \{\tilde{\omega} \dvtx\forall i\in
S, \tilde{\omega}_i = \omega_i\},
\]
and more generally for any $X
\subseteq\Omega$, any $S\dvtx X \to\mathcal{P}(\{1,\ldots,n\})$,
\[
[X]_S:= \bigcup_{\omega\in X} [\omega]_{S(\omega)}.
\]

For any two $A, B \subseteq\Omega$, we denote by $A \circ B$ the
disjoint occurrence of $A$ and $B$ (the notation $A\,\square\,B$ is also
often used to denote this event)
\[
A \circ B := \bigl\{ \omega\dvtx
\mbox{for some $S(\omega) \subseteq\{1,\ldots,n\}$, } [\omega]_S
\subseteq A \mbox{ and } [\omega]_{S^c} \subseteq B \bigr\}.
\]
Recall
that Reimer's inequality states that
%
%e2.5 ###
%
\begin{equation}
\PP(A \circ B) \leq\PP(A) \PP(B).
\end{equation}

We also denote by $\bar{\omega}=1-\omega$ the configuration obtained by
``flipping'' every bit of the configuration $\omega\in\Omega$, and if
$X\subseteq\Omega$, we define $\bar{X} := \{\bar{\omega}\dvtx
\omega\in
X\}$. We are now in a position to state the correlation inequality that
will be a key ingredient in the following:
\begin{theorem}[(\cite{reimerR}, Theorem 1.2)]\label{theo3}
For any $A, B \subseteq\Omega$, we have
%
%e2.6 ###
%
\begin{equation}
\label{ineq}
|A \circ B| \leq|A \cap\bar{B}| = |\bar{A} \cap B|.
\end{equation}
\end{theorem}

For the sake of completeness, let us just mention that this inequality
is not stated explicitly in that form in~\cite{reimerR}. It can be
deduced from Theorem 1.2 by applying it to the \textit{flock of
butterflies} $\mathcal{B} = \{(\omega, f(\omega))\dvtx\omega\in
A\circ
B\}$, where $f(\omega)$ coincides with $\omega$ exactly in the
coordinates in $S(\omega)$. Here, $S(\omega) \subseteq\{1,\ldots,n\}$
is the subset of indices associated with $A$ and $B$ by the definition
of disjoint occurrence, that is, so as to satisfy
$[\omega]_{S(\omega)} \subseteq A$ and $[\omega]_{S(\omega)^c}
\subseteq
B$ for all $\omega\in A \circ B$. For this particular $\mathcal{B}$, we
indeed have $\mathrm{Red} (\mathcal{B}) \subseteq A$ and
$\mathrm{Yellow} (\mathcal{B}) \subseteq B$.

Equivalently, it can be obtained from Lemma 4.1 in~\cite{BCRBKR} by
taking $X = A \circ B$, and $S \dvtx X \to\mathcal{P}(\{1,\ldots,n\})$
associated with $A \circ B$ (so that $[X]_S \subseteq A$ and $[X]_{S^c}
\subseteq B$).

%s2.3 ###
\subsection{Statement of the results}

In this paper, we will be interested in the asymptotic behavior of the
probability of the event $A_{j,\sigma}(n_0(j),N)$ as $N \to\infty$ for
a monochromatic $\sigma$, say $\sigma=B \ldots B$, so that
$A_{j,\sigma}$ simply refers to the existence of $j$ disjoint black
arms. Our first result shows that this probability follows a power law,
as in the case of a polychromatic choice of $\sigma$.
\begin{theorem} \label{thm:existence}
For any $j \geq2$, there exists an exponent $\alpha'_j>0$ such that
%
%e2.7 ###
%
\begin{equation}
\PP(A_{j,B\ldots B}(n,N)) = N^{-\alpha'_j+o(1)}
\end{equation}
as $N \to\infty$ for any fixed $n \geq n_0(j)$.
\end{theorem}

These exponents $\alpha'_j$ are known as the \textit{monochromatic arm
exponents}, and it is natural to try to relate them to the previously
mentioned polychromatic exponents~$\alpha_j$.

Consider any $j \geq2$; we start with a few easy remarks. On the one
hand, one can apply Harris's lemma (more often referred to in the
statistical mechanics community as the FKG inequality, which is its
generalization to certain nonproduct measures; we will follow that
convention in what follows): it implies that
\begin{eqnarray*}
\PP(A_{j+1,B\ldots BW}(n_0,N)) & = &\PP\bigl(A_{j,B\ldots B}(n_0,N) \cap
A_{1,W}(n_0,N)\bigr)\\
& \leq &\PP(A_{j,B\ldots B}(n_0,N)) \PP(A_{1,W}(n_0,N)),
\end{eqnarray*}
and by using item \hyperlink{apriori}{1} above, we get that, for some constants
$C, \varepsilon>0$,
%
%e2.8 ###
%
\begin{equation}
\PP(A_{j+1,B\ldots BW}(n_0,N)) \leq C N^{-\varepsilon}
\PP(A_{j,B\ldots B}(n_0,N))
\end{equation}
or, in other words, that $\alpha'_j < \alpha_{j+1}$. On the other hand,
inequality (\ref{ineq}) directly implies that
\begin{eqnarray*}
\PP(A_{j,B\ldots BB}(n_0,N)) & = & \PP\bigl(A_{j-1,B\ldots B}(n_0,N) \circ
A_{1,B}(n_0,N)\bigr)\\
& \leq &\PP\bigl(A_{j-1,B\ldots B}(n_0,N) \cap A_{1,W}(n_0,N)\bigr)\\
& = &\PP(A_{j,B\ldots BW}(n_0,N)),
\end{eqnarray*}
hence $\alpha'_j \geq\alpha_j$. We will actually prove the following,
stronger result.
\begin{theorem}
\label{thm:strict}
For any $j \geq2$, we have
%
%e2.9 ###
%
\begin{equation}
\alpha_j < \alpha'_j < \alpha_{j+1}.
\end{equation}
\end{theorem}

The monochromatic exponents $\alpha'_j$ thus form a family of exponents
different from the polychromatic exponents.

We would like to stress the fact that the case of half-plane exponents
(or more generally, boundary exponents in any planar domain) is
considerably different: indeed, whenever a boundary is present, the
color-exchange trick implies that the probability of observing $j$ arms
of prescribed colors is \textit{exactly} the same for all color
prescriptions, whether mono- or poly-chromatic. In particular there is
no difference between the monochromatic and polychromatic boundary
exponents. (For the reader's peace of mind, he can notice that the
presence of the boundary provides for a canonical choice of a leftmost
arm; the lack of which is precisely the core idea of the proof of our
main result in the whole plane.)\vadjust{\goodbreak}

We will\vspace*{1pt} first prove the inequality $\alpha_j < \alpha'_j$ (which is the
main statement in the above theorem, the other strict inequality being
the simple consequence of the FKG inequality we mentioned earlier),
since its proof only requires combinatorial arguments, and postpone the
proof of the existence of the exponents to the end of the paper.

In order not to refer to the $\alpha'_j$'s, we adopt the following
equivalent formulation of the inequality: what we formally prove is
that, for any $j \geq2$, there exists $\varepsilon> 0$ such that for
any $N$ large enough,
%
%e2.10 ###
%
\begin{equation}\label{eq:star}
\PP(A_{j,B \ldots BB}(n_0,N)) \leq N^{-\varepsilon} \PP(A_{j,B
\ldots
BW}(n_0,N)).
\end{equation}
The proof of this inequality only relies on discrete features such as
self-duality and RSW-type estimates, and hence it is possible that our
proof could be extended to the case of bond percolation on $\mathbb
Z^2$, where the existence of the exponents, which strongly relies on the
knowledge of the scaling limit, is still unproved; however, the
statement of duality in that case is different enough (the color
exchange trick for instance has no exact counterpart) that some of our
arguments do not seem to extend directly.

%s3 ###
\section{The set of winding angles}

%s3.1 ###
\subsection{Strict inequalities between the exponents}

Our proof is based on an energy versus entropy consideration. The
difference between the monochromatic and the polychromatic $j$-arm
exponents can be written in terms of the expected number of ``really
different'' choices of $j$ arms out of a percolation configuration with
$j$ arms: for a polychromatic configuration, this number is equal to
$1$, whereas for a monochromatic configuration, it grows at least like a
positive power of the modulus, and the ratio between these two numbers
behaves exactly like $(N/n)^{\alpha_j - \alpha'_j}$ because, for fixed
disjoint arms $(r_1, \ldots, r_j)$ with respective lengths $(\ell_1,
\ldots, \ell_j)$, the probability that they are present in the
configuration with a prescribed coloring does not depend on that
coloring (it is equal to $2^{-(\ell_1+\cdots+\ell_j)}$).

\begin{figure}

\includegraphics{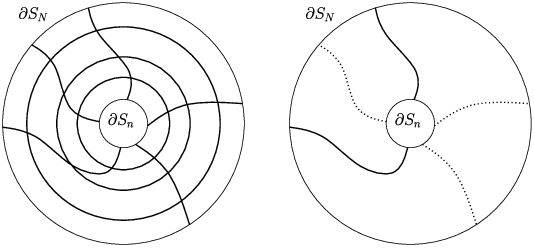}

\caption{To a given monochromatic configuration
correspond many different ``macroscopic'' ways to choose the arms,
contrary to the polychromatic case. The left-hand picture shows the
sort of realization in which many macroscopically different $5$-arm
configurations can be found---the actual topology needed is a
little bit more involved; see Figure \protect\ref{fig:modification2}.}
\label{fig:arms}
\end{figure}

More precisely, but still roughly speaking, the proof relies on the
following observation: given a configuration where $j$ black arms are
present, there are many ways to choose them, since by RSW there is a
positive density of circuits around the origin [allowing ``surgery'' on
the arms (see Figure~\ref{fig:arms})], while if we consider a
configuration with arms of both colors, then there is essentially only
one way to select them. Of course the geometry of an arm is quite
intricate and many local modifications, on every scale, are always
possible, both in the monochromatic and polychromatic setups; what we
mean here is that this choice is unique from a macroscopic point of
view. To formalize this intuition, we thus have to find a way of
distinguishing two macroscopic choices of arms, and for this we will use
the \textit{set of winding angles} associated with a configuration.

%
%f1 ###
%

\begin{definition}
For any configuration of arms, one can choose a continuous
determination of the argument along one of the arms; we call
\textit{winding angle of the arm} (or simply \textit{angle} for
short) the
overall (algebraic) variation of the argument along that arm.
\end{definition}

Clearly, the winding angles of the arms corresponding to a given
$(j,\sigma)$-arm configuration differ by at most $2\pi$. However, for
the same percolation configuration, there might exist many different
choices of a $(j,\sigma)$-arm configuration, corresponding to different
winding angles: we denote by $I_{j,\sigma}(n,N)$ the set of all the
winding angles which can be obtained from such a configuration; we omit
the subscript from the notation whenever $j$ and $\sigma$ are clear from
the context (notice though that whenever there are enough arms for two
different arm events with the same $j$ to be realized simultaneously,
the corresponding sets of winding angles are essentially the same). For
the sake of completeness, we also declare $I_{j,\sigma} (n,N)$ to be
empty if the configuration does not contain $j$ arms of the prescribed
colors.

We will actually rather use $\bar{I}_{j,\sigma}(n,N)$, the set of angles
obtained by ``completing'' $I_{j,\sigma}(n,N)$
\[
\bar{I}_{j,\sigma}(n,N) := \bigcup_{\alpha\in I_{j,\sigma}(n,N)}
( \alpha-\pi,\alpha+\pi].
\]
It is an easy remark that in the
polychromatic case ($\sigma$ nonconstant), we have for any $\alpha\in
I_{j,\sigma}(n,N)$
\[
I_{j,\sigma}(n,N) \subseteq(\alpha-
2\pi,\alpha+2\pi)
\]
(because two arms of different colors cannot cross),
so that $\bar{I}_{j,\sigma}(n,N)$ is an interval of length at most
$4\pi$. In the monochromatic case ($\sigma$ constant), no such bound
applies [and actually it is not obvious that $\bar{I}_{j,\sigma}(n,N)$
is an interval; this is proved as Proposition~\ref{prp:interval}
below].\vadjust{\goodbreak}

In the case of a polychromatic arm configuration, considering successive
annuli of a given modulus as independent, one would expect a central
limit theorem to hold on the angles, or at least fluctuations of order
$\sqrt{\log N}$. On the other hand, for a monochromatic configuration,
performing surgery using circuits in successive annuli should imply that
every time one multiplies the outer radius by a constant, the expected
largest available angle would increase by a constant, so that one would
guess that, by a careful choice of arms, the total angle can be made of
order $\pm\log N$.

Fix $\varepsilon>0$,\vspace*{1pt} and let $A_{1,B}^{\varepsilon}$ (resp.,
$A_{1,W}^\varepsilon$, resp., $A_{j,\sigma}^\varepsilon$) be the event
that there exists a black arm (resp., a white arm, resp., $j$ arms with
colors given by $\sigma$) with angle larger than $\varepsilon\log N$
between radii $n_0$ and $N$. Applying inequality~(\ref{ineq}) with $A =
A_{j-1,B \ldots B}$ and $B = A_{1,B}^\varepsilon$, if the above
intuition was correct, this would imply
\begin{eqnarray*}
\PP(A_{j,B \ldots BB}) & \asymp &\PP(A \circ B)\\
& \leq &\PP(A_{j-1,B \ldots B} \cap A_{1,W}^\varepsilon)\\
& = &\PP(A_{j,B \ldots BW}^\varepsilon),
\end{eqnarray*}
and we could expect
\[
\PP(A_{j,B \ldots BW}^\varepsilon) \leq N^{-\varepsilon'} \PP(A_{j,B
\ldots BW})
\]
by a large-deviation principle. However, proving this LDP seems to be
difficult, and we propose here an alternative proof that relies on the
same ideas, but bypasses some of the difficulties.
\begin{pf*}{Proof of Theorem~\ref{thm:strict}}
We shall in fact prove (\ref{eq:star}). As noted earlier, assuming
Theorem~\ref{thm:existence} (which will be proved in Section
\ref{sec4}), Theorem~\ref{thm:strict} follows immediately.

\textit{Step} 1. First, note that it suffices to prove that
the ratio
\[
\frac{\PP(A_{j,B \ldots BB}(n,N))}{\PP(A_{j,B \ldots
BW}(n,N))}
\]
can be made arbitrarily small as $N/n \to\infty$,
uniformly in $n$: indeed, assuming that this is the case, then for any
$\delta>0$, there exists $R>0$ such that this ratio is less than
$\delta$ as soon as $N/n \geq R$. Then, as a direct consequence of the
quasi-multiplicativity property (item \hyperlink{quasi_mult}{2} above), we
have
\begin{eqnarray*}
&&\PP(A_{j,B \ldots BB}(n,R^k n)) \\
&&\qquad \leq C_2^{k-1} \PP(A_{j,B \ldots BB}(n, R n)) \cdots
\PP(A_{j,B \ldots BB}(R^{k-1}n, R^k n))\\
&&\qquad \leq C_2^{k-1} \delta^k \PP(A_{j,B \ldots BW}(n,
Rn)) \cdots\PP(A_{j,B \ldots BW}(R^{k-1}n, R^{k}n))\\
&&\qquad \leq C_2^{k-1} \delta^k (C_1^{-1})^{k-1} \PP(A_{j,B
\ldots BW}(n, R^{k}n)),
\end{eqnarray*}
and for $\delta= 1/(2 C_2 C_1^{-1})$ this gives
%
%e3.1 ###
%
\begin{equation}
\PP(A_{j,B \ldots BB}(n, R^k n)) \leq2^{-k} \PP(A_{j,B \ldots
BW}(n, R^kn)),\vadjust{\goodbreak}
\end{equation}
which immediately implies that for some $C, \varepsilon> 0$,
\[
\PP(A_{j,B \ldots BB}(n,N)) \leq C \biggl( \frac{N}{n} \biggr)
^{-\varepsilon} \PP(A_{j,B \ldots BW}(n,N)).
\]
In particular, applying
this for $n=n_0$ (and $N$ large enough) leads to the inequality that
we need.

\textit{Step} 2. The key step of the proof is as follows.
Given a configuration with $j$ arms in an annulus of large modulus, we
use RSW-type estimates to prove the existence of a large number of
disjoint sub-annuli of it, in each of which one can find black paths
topologically equivalent to those in Figure~\ref{fig:modification} (in
the case $j=2$) or its reflection. Every time this configuration
appears, one has the possibility to replace the original arms (in
solid lines on the figure) with modified, and still disjoint, arms,
obtained by using one of the dashed spirals in each of them. The
new arms then land at the same points on the outer circle, but with a
winding angle differing by $2\pi$. This allows us to show that, with
high probability, the set of angles $\bar{I}(n,N)$ contains an
interval of length at least $\varepsilon\log(N/n)$, for some
$\varepsilon>0$ (which can be written in terms of the RSW estimates). We
now proceed to make the construction in detail.

%
%f2 ###
%
\begin{figure}

\includegraphics{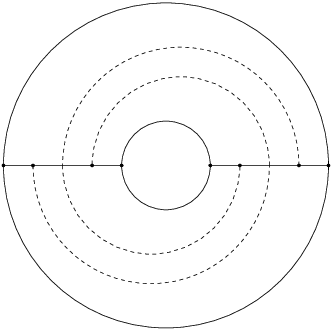}

\caption{When they encounter this
configuration, the arms (here in solid lines) can be modified,
detouring via the dashed lines, to make an extra turn.}\label{fig:modification}
\end{figure}

%
%f3 ###
%
\begin{figure}

\includegraphics{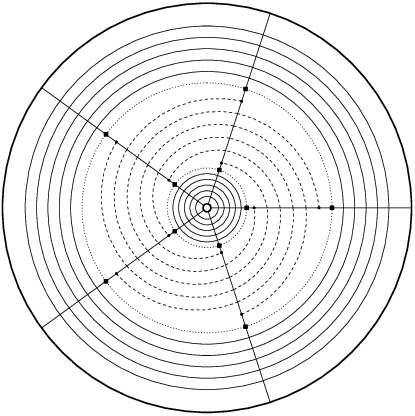}

\caption{Generalization of
Figure \protect\ref{fig:modification} in the case of $j \geq3$ arms. The
additional circuits (in solid lines) are needed to apply Menger's
theorem; the circles of radii $m$ and $4m$ (resp., $2m$ and $3m$)
are drawn in heavy (resp., dotted) lines, the spiraling paths in
dashed lines and the active points are marked with a black
square.}\label{fig:modification2}\vspace*{-2pt}
\end{figure}

Let $j\geq2$, and let $m$ be a positive integer. Define a
\textit{$j$-spiral} between radii $m$ and $4m$ as the configuration
pictured in Figure~\ref{fig:modification2}. More precisely, a
$j$-spiral is the union of $4$ families of $j$ black paths in a
percolation configuration, namely:
\begin{itemize}
\item$j$ disjoint rays between radii $m$ and $4m$;
\item$j$ disjoint ``spiraling paths'' contained in the annulus
$S_{2m,3m}$, each connecting two points of one of the rays and
making one additional turn around the origin;\vadjust{\goodbreak}
\item$j$ disjoint circuits around the origin, contained in the
annulus $S_{m,2m}$;
\item$j$ disjoint circuits around the origin, contained in the
annulus $S_{3m,4m}$.
\end{itemize}
RSW-type estimates directly show that, uniformly as $m\to\infty$, the
probability of observing a $j$-spiral between radii $m$ and $4m$ is
bounded below by a positive constant (depending only on $j$). In
addition, with each such spiral we associate two families of $j$
\textit{active points}: for each of the $j$ rays, oriented starting at
radius $m$, its last intersection with the circle of radius $2m$ is
called an \textit{inner} active point, and its first intersection with
the circle of radius $3m$ after it, an \textit{outer} active point; in
particular, the section of the ray between its two active points
remains within the annulus $S_{2m,3m}$.

The same RSW arguments show that with positive probability, such a
$j$-spiral actually satisfies a few additional properties (which we
will consider part of the definition from now on): the different
pieces remain well separated whenever they are not forced to intersect
by topological constraints; whenever two pieces intersect, they can do
so multiple times, but all intersection points remain close to each
other, and it is easy to check that this is not a problem in the proof.
In fact, the only thing that might cause trouble is that the rays may
cross each of the circles of radii $2m$ and $3m$ several times; one
cannot use RSW to avoid it, but it is taken care of by our choice of
the active points.

The presence of $j$-spirals in disjoint annuli are independent events,
each with positive probability,\vadjust{\goodbreak} so that,\vspace*{1pt} for some $\varepsilon>0$, the
probability of the event $E^{(\varepsilon)}_j(n,N)$ of having at least
$\varepsilon\log(N/n) + 1$ disjoint $j$-spirals between radii $n$ and
$N$ goes to~$1$ as $N/n$ goes to infinity. The presence of $j$-spirals
being an increasing event, the FKG inequality ensures the conditional
probability of $E^{(\varepsilon)}_j(n,N)$, given the existence of $j$
black arms between radii $n$ and $N$, still goes to $1$ as $N/n$ goes
to infinity.

We now explain how to use $j$-spirals to perform surgery on black
arms. Assume that between radii $n$ and $N$, there are a certain
number $s \geq2$ of disjoint spiral configurations; let $(m_i)_{1
\leq i \leq s}$ be the corresponding values of $m$. Assume in addition
that there are $j$ disjoint arms between radii $n$ and $N$. The first
remark is the following: for every $i \in\{1, \ldots, s-1\}$, the
event $A_{j,B\ldots BB}(3m_i, 2m_{i+1})$ is realized, and one can
choose $j$ disjoint arms $(\gamma_{i,k})_{1\leq k\leq j}$ in the
annulus $S_{3m_i, 2m_{i+1}}$ accordingly. The main part of the
argument then consists of proving that, within the union of all those
arms together with the spirals, it is always possible to find
$2^{s-1}$ $j$-arm configurations between radii $n$ and $N$, spanning
winding angles in an interval of length $2(s-2)\pi$.

Each of the constructed arms will consist of sections of two types,
connecting at the active points (and disjoint otherwise, by our choice
of the active points): the ``movable'' ones, formed out of the rays
and spiraling parts of the $j$-spirals and contained in the union of
the annuli $S_{2m_i,3m_i}$; and the ``intermediate'' ones, made from
the arms $\gamma_{i,k}$ (which we have chosen to stay in the annuli
$S_{3m_i,2m_{i+1}}$) and from pieces of the $j$-spirals lying outside
the active points. Within every $j$-spiral there are two choices for
the corresponding movable section, which is how arms of very different
winding angles will be produced. For that, all we need to prove is the
existence of one family of $j$ disjoint arms following this
decomposition.

The key argument is as follows. Define $\Gamma_1$ as the union of the
$(\gamma_{1,k})_{1 \leq k \leq j}$ (which cross the annulus
$S_{3m_1, 2m_2}$), together with the parts of the $j$-spiral
$\Sigma_1$ (resp., $\Sigma_2$) outside its outer active points (resp.,
inside its inner active points). It is easy to check that, whenever
one marks $(j-1)$ vertices on $\Gamma_1$, there still exists a path
completely contained in $\Gamma_1$ and avoiding the marked points, and
connecting one of the outer active points of $\Sigma_1$ to one of the
inner active points of $\Sigma_2$. Indeed, marking $(j-1)$ vertices
leaves untouched at least one arm, one inner (resp., outer) circle
and one ray of $\Sigma_1$ (resp., $\Sigma_2$). Menger's theorem (see
\cite{diestelbook}) then ensures that $\Gamma_1$ contains $j$
disjoint arms, each connecting one of the outer active points of
$\Sigma_1$ to one of the inner active points of $\Sigma_2$.

The same construction can be performed between each pair of successive
$j$-spirals, to produce the desired ``intermediate'' sections. An
obvious variation of the construction can also be applied in the
annuli $S_{n,2m_1}$ and $S_{3m_s,N}$, leading to the following fact:
whenever there are $j$ arms between radii $n$ and $N$, and the event
$E_j^{(\varepsilon)}(n,N)$ is realized, the set $\bar I(n,N)$ contains an
interval of length at least $2\pi\varepsilon\log(N/n)$ (obtained by
playing with the appropriate movable sections), and this occurs
with conditional probability going to $1$ as $N/n$ goes to infinity.
Notice that the winding angles of the original $j$ arms need not be
within the interval we just constructed; however, this has no bearing
on what follows.

\textit{Step} 3. We now use the BK inequality to control the
probability, given the presence of $j$ black arms between radii $n$
and $N$, that there is a choice of arms with a very large winding
angle. In fact, the argument shows a little more: it is very unlikely,
even without the conditioning, that there is a single arm with large
winding.

Let $R_{m,m'}$ be the intersection of the annulus $S_{m,m'}$ with the
cone $\mathcal C := \{z \in\mathbb C \dvtx|{\arg(z)}| <
\pi/10\}$ (using the standard identification of the plane $\mathbb
R^2$ with the complex plane). We will consider two families of
``rectangles'': the $R_{e^k, e^{k+2}}$ (which we call the \textit{long}
ones), and the $R_{e^k, e^{k+1}}$ (the \textit{wide} ones). It is easy
to see that any curve connecting two points of the plane of arguments
$- \pi/10$ and $+ \pi/10$ while staying within the cone $\mathcal C$
has to cross at least one of these rectangles between two opposite
sides. More precisely, it must either cross a long rectangle from one
straight side to the other, or a wide rectangle from one curved side
to the other.

Now, assume that there exist $j$ arms between radii $n$ and $N$, and
that their winding angle is at least equal to $2 \pi K \log(N/n) +
4\pi$ (where $K$ is a positive constant which will be chosen later).
Each of these arms has to cross the cone $\mathcal C$ at least (the
integer part of) $K \log(N/n)$ times, so the configuration contains
$jK \log(N/n)$ disjoint paths (at least), each one crossing one of the
rectangles. On the other hand, the probability that a rectangle of a
given shape is crossed by a path is bounded above by $1-\delta$ for
some $\delta>0$, as provided by RSW estimates.

Combined with the BK inequality, this implies that the probability
that there are $j$ arms winding of at least an angle of $2 \pi K \log
(N/n) + 4\pi$ is bounded above by
\[
\mathcal S := \sum_{(l_k,l'_k)}
\prod_k (1-\delta)^{l_k + l'_k},
\]
where the sum is taken over all
$\log(N/n)$-tuples of $(l_k,l'_k)$ having a sum equal to $jK
\log(N/n)$. The number of such tuples is the same as the number of
choices of $2\log(N/n)-1$ disjoint elements out of
$(jK+2)\log(N/n)-1$, so we obtain that
\[
\mathcal S \leq
(1-\delta)^{jK \log(N/n)} \pmatrix{(jK+2) \log(N/n)-1\cr
2\log(N/n)-1}.
\]
It is then a straightforward application of
Stirling's formula to obtain that
\[
\mathcal S \leq C \exp
[ (-c K + C \log K) \log(N/n) ],
\]
where the constants $c$ and
$C$ do not depend on the value of $K$. Choosing $K$ large enough, one
then obtains that
\[
\mathcal S \leq C \biggl( \frac n N
\biggr)^{\beta_j+1}
\]
[where $\beta_j$ is the same as in
(\ref{apriorij})].

Hence, the conditional probability, given $A_{j,B \ldots BB}(n,N)$,
that $\bar I(n,N)$ is contained in the interval of length $4\pi K \log
(N/n)+8\pi$ centered at $0$ goes to $1$ as $N/n$ goes to infinity.
Dividing that interval into sub-intervals of length
$\frac{\varepsilon}{2} \log(N/n)$, and using the previous step, we get
that for one of them, say $i_{\varepsilon}(n,N)$,
\[
\PP\bigl(i_{\varepsilon}(n,N)
\subseteq\bar{I}(n,N) | A_{j,B \ldots BB}(n,N)\bigr) \geq C',
\]
where
$C'>0$ is a universal constant.

\textit{Step} 4. We are now in a position to conclude. If we
take $\alpha_{\min}$ such that
\[
\PP\bigl( A_{j,B \ldots
BW}(n,N) \cap\{ \alpha_{\min} \in\bar{I}_{j,B\ldots BW}(n,N)
\} \bigr)
\]
is minimal among all $\alpha_{\min} \in
i_{\varepsilon}(n,N) \cap(4 \pi\mathbb{Z})$, then
\begin{eqnarray*}
&&\PP\bigl(A_{j,B \ldots BW}(n,N) \cap \{ \alpha_{\min} \in
\bar{I}_{j,B\ldots BW}(n,N) \}\bigr) \\
&&\qquad \leq
\frac{4\pi}{{\varepsilon}/{2} \log(N/n)} \PP(A_{j,B \ldots
BW}(n,N))
\end{eqnarray*}
since, as we noted earlier, whenever there are arms of different
colors, $4 \pi\mathbb{Z}$ cannot contain more than one element of
$\bar{I}(n,N)$. On the other hand, we know from the previous step that
\begin{eqnarray*}
&&
\PP\bigl(A_{j,B \ldots BB}(n,N) \cap\{ \alpha_{\min} \in
\bar{I}_{j,B\ldots BB}(n,N) \}\bigr) \\
&&\qquad \geq\PP\bigl(A_{j,B \ldots BB}(n,N)
\cap\{ i_{\varepsilon}(n,N) \subseteq\bar{I}_{j,B\ldots BB}(n,N)
\}\bigr)\\
&&\qquad \geq C' \PP(A_{j,B \ldots BB}(n,N)).
\end{eqnarray*}
If we apply (\ref{ineq}) to $A = A_{j-1,B \ldots B}(n,N)
\cap\{ \alpha_{\min} \in\bar{I}_{j-1, B\ldots B}(n,N) \}$ and
$B = A_{1,B}(n,N)$, we obtain that
\begin{eqnarray*}
C' \PP(A_{j,B \ldots BB}(n,N)) &\leq& \PP\bigl(A_{j,B \ldots
BB}(n,N) \cap\{ \alpha_{\min} \in\bar{I}_{j,B\ldots B}(n,N)
\}\bigr)\\
& = & \PP(A \circ B)\\
& \leq& \PP(A \cap\bar{B})\\
& = & \PP\bigl(A_{j,B \ldots BW}(n,N) \cap\{ \alpha_{\min} \in
\bar{I}_{j-1,B\ldots B}(n,N) \}\bigr)\\
&\leq& \PP\bigl(A_{j,B \ldots BW}(n,N) \cap\{ \alpha_{\min} \in
\bar{I}_{j,B\ldots BW}(n,N) \}\bigr)\\
&\leq& \frac{4\pi}{{\varepsilon}/{2} \log(N/n)} \PP(A_{j,B
\ldots BW}(n,N)),
\end{eqnarray*}
which completes the proof.
\end{pf*}

%s3.2 ###
\subsection{The density of the set of winding angles}

In this section, we further describe the set of winding angles
$I(n,N)$, which happened to be a key tool in the previous proof, in the
monochromatic case. We prove that (conditionally on the existence of $j$
disjoint black arms) $\bar{I}(n,N)$ is always an interval, as in the
polychromatic case. For that, we use the following deterministic
statement that $I(n,N)$ does not have large ``holes'':
\begin{proposition} \label{prp:interval}
Let $j \geq1$ and $\sigma= B \ldots BB$ of length $j$. Let
$\alpha,\alpha'\in I_{j,\sigma}(n,N)$ with $\alpha< \alpha'$; then
there exists a sequence $(\alpha_i)_{0\leq i\leq r}$ of
elements of $I_{j,\sigma}(n,N)$, satisfying the following two
properties:
\begin{itemize}
\item$\alpha= \alpha_0 < \alpha_1 < \cdots< \alpha_r = \alpha'$;
\item for every $i \in\{0,\ldots,r-1\}$, $\alpha_{i+1} - \alpha_i <
2\pi$.
\end{itemize}
\end{proposition}

This result directly implies that $\bar{I}(n,N)$ is an interval, and the
construction of the previous subsection, creating extra turns (step 2
of the proof), gives a lower bound on the diameter of $\bar{I}(n,N)$: we
hence get that for $\sigma$ constant, there exists some $\varepsilon>0$
(depending only on $j$) such that $\bar{I}(n,N)$ is an interval of
length at least $\varepsilon\log(N/n)$ with probability tending to
$1$ as
$N/n$ gets large.

The main step in the proof of the density result is the following
topological lemma:
\begin{lemma}\label{lemmaa}
Let $j \geq1$, and let $\gamma_1, \ldots, \gamma_j$ be $j$ disjoint
Jordan curves contained in the (closed) annulus $\{n \leq|z| \leq
N\}$, ordered cyclically and each having its starting point on the
circle of radius $n$ and its endpoint on the circle of radius $N$. For
each $k \in\{1,\ldots,j\}$, let $\alpha_k$ be the winding angle of
$\gamma_k$ (as defined above) and let $\delta_k$ be the ray $[n e^{2 i
\pi k/j}, N e^{2 i \pi k/j}]$. Assume that, for each pair $(k,k')$,
the intersection of $\gamma_k$ and $\delta_{k'}$ is finite. Then,
provided all the $\alpha_k$ are larger than $2\pi(1+2/j)$, the union
of all the paths $\gamma_k$ and $\delta_k$ contains $j$ disjoint paths
$\tilde\delta_1, \ldots, \tilde\delta_j$, all having winding angle
$2\pi/j$ and sharing the same endpoints as the $\delta_k$.
\end{lemma}

In other words: starting from two collections of paths, if their angles
differ enough, one can ``correct'' the one with the smaller angle in
such a way as to make it turn a little bit more.
\begin{pf*}{Proof of Lemma~\ref{lemmaa}}
We shall construct the paths $\tilde\delta_k$ explicitly. The first
step is to reduce the situation to one of lower combinatorial
complexity, namely to the case where the starting points of the
$\gamma_k$ are separated by those of the $\delta_k$. For each
$k\leq j$, let $\tau_k = \inf\{t \dvtx\gamma_k(t) \in[n e^{i \pi
(2k-1)/j}, N e^{i \pi(2k-1)/j}]$\} (which is always finite by our
hypotheses), and let
\[
\Gamma:= \bigcup_{k=1}^j \{\gamma_k(t) \dvtx0
\leq t \leq\tau_k \}.
\]
$\Gamma$ intersects each of the
$\delta_k$ finitely many times, so each of the $\delta_k \setminus
\Gamma$ has finitely many connected components: let $\Delta$ be the
union of those components that do not intersect the circle of radius
$N$, and let
\[
\Omega_0 := \{ n \leq|z| \leq N\} \setminus(\Gamma
\cup\Delta).
\]

Let $\Omega$ be the connected component of $\Omega_0$ having the
circle of radius $N$ as a boundary component. $\Omega$ is homeomorphic
to an annulus, and for each $k$, the point $\gamma_k(\tau_k)$ is on
its boundary; by construction, the $\gamma_k(\tau_k)$ are intertwined
with the (remaining portions of the) rays of angles $2\pi k/j$. We
will perform our construction of the $\tilde\delta_k$ inside
$\Omega$; continuing them with the $\delta_k$ outside $\Omega$ then
produces $j$ disjoint paths satisfying the conditions we need.

Up to homeomorphism, we can now assume without loss of generality that
for each $k$, $\gamma_k(0) = n e^{i \pi(2k-1)/j}$. The only thing we
lose in the above reduction is the assumption on the angles of the
$\gamma_k$; but since it takes at most one turn for each of the
$\gamma_k$ to reach the appropriate argument, we can still assume that
the remaining angles are all larger than $4\pi/j$. In particular, each
of the $\gamma_k$ will cross the wedge between angles $2\pi k/j$ and
$2\pi(k+1)/j$ in the positive direction before hitting the circle of
radius $N$.

For every $k \leq j$, let $\theta_k(t)$ be the continuous
determination of the argument of $\gamma_k(t)$ satisfying $\theta_k(0)
= (2k-1)\pi/j$, and let
\[
\mathcal T_k := \biggl\{ t>0 \dvtx\frac{2\pi
k}j < \theta_k(t) < \frac{2\pi(k+1)}j \biggr\}\quad \mbox{and}\quad
\tilde\Gamma_k = \overline{ \{ \gamma_k(t) \dvtx t \in\mathcal
T_k\} }.
\]

We now describe informally the construction of $\tilde\delta_k$. Start
from the point $n e^{2 \pi i k/j}$, and start following $\delta_k$
outward, until the first intersection of $\delta_k$
with~$\tilde\Gamma_k$. Then, follow the corresponding connected
component of $\tilde\Gamma_k$, until intersecting either $\delta_k$ or
$\delta_{k+1}$; follow that one outward until it intersects either
$\tilde\Gamma_k$ or the circle of radius $N$; iterating the
construction, one finally obtains a Jordan path\vspace*{1pt} joining $n e^{2\pi i
k/j}$ to $N e^{2\pi i (k+1)/j}$, and contained in the union of
$\delta_k$, $\delta_{k+1}$ and $\tilde\Gamma_k$ (see
Figure~\ref{fig:deltatilde}).

%
%f4 ###
%
\begin{figure}

\includegraphics{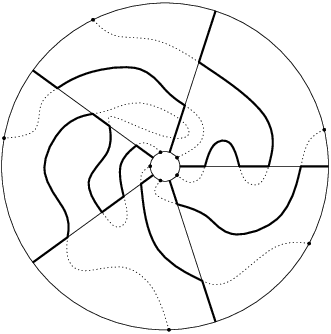}

\caption{The construction of the $\tilde
\delta_k$ (in the case $j=5$). The dotted lines are the paths
$\gamma_k$, and the heavy lines are the $\tilde\delta_k$ obtained
at the end of the construction.}\label{fig:deltatilde}
\end{figure}

All that remains is to prove that the $\tilde\delta_k$ are indeed
disjoint; by symmetry, it is enough to do so for $\tilde\delta_1$ and
$\tilde\delta_2$. Besides,\vspace*{1pt} because the $\gamma_k$ are themselves
disjoint, any intersection point between $\tilde\delta_1$ and $\tilde
\delta_2$ has to occur on $\delta_2$ (at least in the case $j>2$, but
the case $j=2$, where they could also intersect along $\delta_1$,
again follows by symmetry).

The intersection of $\tilde\delta_1$ with $\delta_2$ consists of a
finite collection $(I_m)$ of compact intervals; besides, the points of
the intersection are visited by $\tilde\delta_1$ in order of
increasing distance to the origin. Similarly, the intersection of
$\tilde\delta_2$ with $\delta_2$ consists of a finite collection
$(J_l)$ of compact intervals, which are also visited in order of
increasing distance to the origin.

Suppose that $\bigcup I_p$ and $\bigcup J_p$ have a nonempty
intersection, and let $z_0$ be the intersection point lying closest to
the origin. Let $p_0$ and $q_0$ be such that $z_0 \in I_{p_0} \cap
J_{q_0}$; notice that $z_0$ is the endpoint closest to the origin of
either $I_{p_0}$ or $J_{q_0}$. According to the order in which
$\gamma_1$ (resp., $\gamma_2$) visits the endpoints of $I_{p_0}$
(resp., $J_{q_0}$), this gives rise to eight possible configurations;
it is straightforward in all cases to apply Jordan's theorem to prove
that $\gamma_1$ and $\gamma_2$ then have to intersect, thus leading to
a contradiction.
\end{pf*}

For the purpose of the proof of Proposition~\ref{prp:interval}, we will
need a slight variation of the lemma, where the hypothesis of finiteness
of the intersections between paths is replaced with the assumption that
the paths considered are all polygonal lines. The proof is exactly the
same though, and does not even require any additional notation: whenever
two paths, say $\gamma_k$ and $\delta_{k'}$, coincide along a line
segment, the definition of $\Gamma_k$ amounts to considering some of the
endpoints of this segment as intersections, which in other words is
equivalent to shifting $\gamma_k$ by an infinitesimal amount toward the
exterior of the wedge used to define $\mathcal T_k$ in order to recover
finiteness.
\begin{pf*}{Proof of Proposition~\ref{prp:interval}}
The previous lemma is stated with particular curves on which a surgery
can be done, but it can obviously be applied to more general cases
through a homeomorphism of the annulus. The general statement is then
the following (roughly speaking): assuming the existence of two
families of $j$ arms with different enough winding angles, it is
possible to produce a third family using the same endpoints as the
first one but with a slightly larger winding angle.

We are now ready to prove Proposition~\ref{prp:interval}. Consider a
configuration in which one can find two families of crossings, say
$(\lambda_k)$ and $(\lambda'_k)$, in such a way that for every~$k$,
the difference between the winding angles of $\lambda_k$ and
$\lambda'_k$ is at least $2\pi$. Let $\alpha_0$ be the minimal angle
in the first family, and apply the topological lemma with $\delta_k =
\lambda_k$ and $\gamma_k = \lambda'_k$: one obtains a new family of
pairwise disjoint paths $(\lambda^1_k)$, which share the same family
of endpoints as the $(\lambda_k)$, the endpoint of $\lambda_k^1$ being
that of $\lambda_{k+1}$ (with the obvious convention that $j+1=1$).

One can then iterate the procedure, applying the topological lemma
with this time $\delta_k = \lambda^1_k$, and still letting $\gamma_k =
\lambda'_k$; one gets a new family $(\lambda^2_k)$ with the endpoints
again shifted amongst the paths in the same direction. Continuing as
long as the winding angle difference is at least $2\pi$, this
construction produces a sequence $(\lambda^i_k)$ of $j$-tuples of
disjoint paths, the winding angles of which vary by less than $2\pi$
at each step. Besides, the construction ends in finitely many steps,
for after $j$ steps, each of the winding angles has increased by
exactly $2\pi$. This readily implies our claim.
\end{pf*}
\begin{remark}
Notice that, as early as the second step of the procedure,
$(\lambda^n_k)$~and $(\lambda'_k)$ will always coincide on a positive
fraction of their length, which is why we needed the above extension
of the lemma.
\end{remark}

%s4 ###
\section{Existence of the monochromatic exponents}\label{sec4}

We now prove Theorem~\ref{thm:existence}, stating the existence of the
monochromatic exponents $\alpha'_j$. For that, we use a rather common
argument, presented, for example, in~\cite{smirnovexps}: since the
quasi-multiplicativity property holds (item \hyperlink{quasi_mult}{2} above),
it is actually enough to check that there exists a function $f_j$ (which
will automatically be sub-multiplicative itself; one can take $C_2=1$
in the quasi-multiplicativity property) such that, for every $R>1$,
%
%e4.1 ###
%
\begin{equation}
\label{eq:convergence}
\PP(A_{j,B \ldots BB}(n,R n)) \to f_j(R)
\end{equation}
as $n \to\infty$. Notice that RSW-type estimates provide both the fact
that the left-hand term in bounded above and below by constants for
fixed $R$ as $n\to\infty$, and a priori estimates on any
(potentially subsequential) limit, of the form
\[
R^{-1/\varepsilon_j} \leq
f_j(R) \leq R^{-\varepsilon_j},
\]
where $\varepsilon_j$ depends only on $j$.

By Menger's theorem (see, e.g.,~\cite{diestelbook}), the
complement of the event $A_{j,B \ldots BB}(n$, $N)$ can be written as
\begin{eqnarray*}
D_{j}(n,N) &=& \{ \mbox{There exists a circuit in $S_{n,N}$
that surrounds $\partial^i S_n$} \\
&&\hspace*{52.6pt} \mbox{and contains at most $j-1$ black sites} \}.
\end{eqnarray*}
This makes it possible to express the event $A_{j,B \ldots BB}(n,N)$ in
terms of the collection of all cluster interfaces (or ``loops''): it is
just the event that there does not exist a ``necklace'' of at most
$(j-1)$ loops, with white vertices on their inner boundary and black
ones on their outer boundary, forming a chain around $\partial^i S_n$
and such that two consecutive loops are separated by only one black
site.

Standard arguments show that the probability that two interfaces touch
in the scaling limit is exactly the asymptotic probability that they
``almost touch'' (in the sense that they are separated by exactly one
vertex) on discrete lattices; it is, for example, a simple consequence
of the fact that the polychromatic six-arm exponent is strictly larger
than $2$, which in turn is a consequence of RSW-type estimates (the fact
that the polychromatic five-arm exponent is equal to $2$ being true on
any lattice on which RSW holds, at least for colors $BWBWW$).

What this means, is that to show convergence of the probability in
(\ref{eq:convergence}), it is enough to know the probability of
the corresponding continuous event. While we do not know the exact value
of the limit, it is nevertheless easy to check that the event itself is
measurable with respect to the \textit{full scaling limit} of percolation,
as constructed by Camia and Newman in~\cite{camiafull}, and that is
enough for our purpose. Notice that the measurability of the event in
terms of the full scaling limit is ensured by the exploration procedure
described in that paper: it is proved there that for every
$\varepsilon>0$, all loops of diameter at least $\varepsilon$ (with the
proper orientation) are discovered after finitely many steps of the
exploration procedure.
\begin{remark}
Proving the existence of the exponents only requires the existence of
the function $f_j$, but it is easily seen that in fact the value of
$\alpha'_j$ describes the power law decay of $f_j(R)$ as $R$ goes to
infinity. However, deriving the value of the exponent directly from
the full scaling limit seems to be difficult, and we were not able to
do it.
\end{remark}

\section*{Acknowledgments}

We are indebted to J. van den Berg for pointing out inequality
(\ref{ineq}) to us. We are also very grateful to W. Werner for many
stimulating discussions, and to an anonymous referee for many
insightful remarks about the first version of this paper. Part of this
research was realized during a semester (resp., year) spent by Pierre
Nolin (resp., Vincent Beffara) at Universit\'{e} de Gen\`{e}ve, and
both authors would like to thank the mathematics department there for
its hospitality, in particular Stanislav Smirnov.

\printaddresses


\begin{thebibliography}{99}

%b1 ###
\bibitem{aharonyexponents}
%
\begin{barticle}[author]
\bauthor{\bsnm{Aizenman},~\bfnm{M.}\binits{M.}},
\bauthor{\bsnm{Duplantier},~\bfnm{Bertrand}\binits{B.}} \AND
\bauthor{\bsnm{Aharony},~\bfnm{A.}\binits{A.}}
(\byear{1999}).
\btitle{Path crossing exponents and the external perimeter in {2D}
percolation}.
\bjournal{Phys. Rev. Lett.}
\bvolume{83}
\bpages{1359--1362}.
\end{barticle}
%
\endbibitem

%b2 ###
\bibitem{BCRBKR}
%
\begin{binproceedings}[vtex]
\bauthor{\bsnm{Borgs},~\bfnm{Christian}\binits{C.}},
\bauthor{\bsnm{Chayes},~\bfnm{Jennifer}\binits{J.}} \AND
\bauthor{\bsnm{Randall},~\bfnm{Dana}\binits{D.}}
(\byear{1999}).
\btitle{The van den Berg--Kesten--Reimer inequality: A review}.
In \bbooktitle{Perpexing Problems in Probability: Festschrifft in
honor of Harry Kesten}.
\bpublisher{Birkh\"auser}, \baddress{Boston, MA}.
\bid{mr={1703130}}
\end{binproceedings}
%
\endbibitem

%b3 ###
\bibitem{camiafull}
%
\begin{barticle}[author]
\bauthor{\bsnm{Camia},~\bfnm{Federico}\binits{F.}} \AND
\bauthor{\bsnm{Newman},~\bfnm{Charles~M.}\binits{C.~M.}}
(\byear{2005}).
\btitle{Two-dimensional critical percolation: The full scaling limit}.
\bjournal{Comm. Math. Phys.}
\bvolume{268}
\bpages{1--38}.
\bid{doi={10.1007/s00220-006-0086-1}, mr={2249794}}
\end{barticle}
%
\endbibitem

%b4 ###
\bibitem{diestelbook}
%
\begin{bbook}[vtex]
\bauthor{\bsnm{Diestel},~\bfnm{R.}\binits{R.}}
(\byear{2000}).
\btitle{Graph Theory}, \bedition{2nd} ed.
\bpublisher{Springer}, \baddress{New York}.
\bid{doi={10.1007/b100033}, mr={1743598}}
\end{bbook}
%
\endbibitem

%b5 ###
\bibitem{grimmettbook}
%
\begin{bbook}[vtex]
\bauthor{\bsnm{Grimmett},~\bfnm{Geoffrey~R.}\binits{G.~R.}}
(\byear{1999}).
\btitle{Percolation},
\bedition{2nd} ed.
\bseries{Grundlehren der Mathematischen Wissenschaften [Fundamental Principles
  of Mathematical Sciences]}
\bvolume{321}.
\bpublisher{Springer}, \baddress{Berlin}.
\bid{mr={1707339}}
\end{bbook}
%
\endbibitem

%b6 ###
\bibitem{kestenscaling}
%
\begin{barticle}[vtex]
\bauthor{\bsnm{Kesten},~\bfnm{Harry}\binits{H.}}
(\byear{1987}).
\btitle{Scaling relations for {2D}-percolation}.
\bjournal{Comm. Math. Phys.}
\bvolume{109}
\bpages{109--156}.
\bid{mr={879034}}
\end{barticle}
%
\endbibitem


%b7 ###
\bibitem{wernervalue}
%
\begin{barticle}[vtex]
\bauthor{\bsnm{Lawler},~\bfnm{Gregory~F.}\binits{G.~F.}},
\bauthor{\bsnm{Schramm},~\bfnm{Oded}\binits{O.}} \AND
\bauthor{\bsnm{Werner},~\bfnm{Wendelin}\binits{W.}}
(\byear{2001}).
\btitle{Values of {Brownian} intersection exponents {I}: Half-plane exponents}.
\bjournal{Acta Math.}
\bvolume{187}
\bpages{237--273}.
\bid{doi={10.1007/BF02392618}, mr={1879850}}
\end{barticle}
%
\endbibitem

%b8 ###
\bibitem{wernervalue2}
%
\begin{barticle}[vtex]
\bauthor{\bsnm{Lawler},~\bfnm{Gregory~F.}\binits{G.~F.}},
\bauthor{\bsnm{Schramm},~\bfnm{Oded}\binits{O.}} \AND
\bauthor{\bsnm{Werner},~\bfnm{Wendelin}\binits{W.}}
(\byear{2001}).
\btitle{Values of {Brownian} intersection exponents {II}: Plane exponents}.
\bjournal{Acta Math.}
\bvolume{187}
\bpages{275--308}.
\bid{doi={10.1007/BF02392619}, mr={1879851}}
\end{barticle}
%
\endbibitem

%b9 ###
\bibitem{lswonearm}
%
\begin{barticle}[vtex]
\bauthor{\bsnm{Lawler},~\bfnm{Gregory~F.}\binits{G.~F.}},
\bauthor{\bsnm{Schramm},~\bfnm{Oded}\binits{O.}} \AND
\bauthor{\bsnm{Werner},~\bfnm{Wendelin}\binits{W.}}
(\byear{2002}).
\btitle{One-arm exponent for critical 2{D} percolation}.
\bjournal{Electron. J. Probab.}
\bvolume{7}
\bpages{13 pp. (electronic)}.
\bid{mr={1887622}}
\end{barticle}
%
\endbibitem

%b10 ###
\bibitem{reimerR}
%
\begin{barticle}[vtex]
\bauthor{\bsnm{Reimer},~\bfnm{D.}\binits{D.}}
(\byear{2000}).
\btitle{Proof of the van den Berg--Kesten conjecture}.
\bjournal{Combin. Probab. Comput.}
\bvolume{9}
\bpages{27--32}.
\bid{doi={10.1017/S0963548399004113}, mr={1751301}}
\end{barticle}
%
\endbibitem

%b11 ###
\bibitem{schrammlerw}
%
\begin{barticle}[vtex]
\bauthor{\bsnm{Schramm},~\bfnm{Oded}\binits{O.}}
(\byear{2000}).
\btitle{Scaling limits of loop-erased random walks and uniform
spanning trees}.
\bjournal{Israel J. Math.}
\bvolume{118}
\bpages{221--288}.
\bid{doi={10.1007/BF02803524}, mr={1776084}}
\end{barticle}
%
\endbibitem


%b12 ###
\bibitem{smirnovperco}
%
\begin{barticle}[vtex]
\bauthor{\bsnm{Smirnov},~\bfnm{Stanislav}\binits{S.}}
(\byear{2001}).
\btitle{Critical percolation in the plane: Conformal invariance, {Cardy}'s
formula, scaling limits}.
\bjournal{C. R. Acad. Sci. Paris S\'er. I Math.}
\bvolume{333}
\bpages{239--244}.
\bid{doi={10.1016/S0764-4442(01)01991-7}, mr={1851632}}
\end{barticle}
%
\endbibitem

%b13 ###
\bibitem{smirnovexps}
%
\begin{barticle}[vtex]
\bauthor{\bsnm{Smirnov},~\bfnm{Stanislav}\binits{S.}} \AND
\bauthor{\bsnm{Werner},~\bfnm{Wendelin}\binits{W.}}
(\byear{2001}).
\btitle{Critical exponents for two-dimensional percolation}.
\bjournal{Math. Res. Lett.}
\bvolume{8}
\bpages{729--744}.
\bid{mr={1879816}}
\end{barticle}
%
\endbibitem

%b14 ###
\bibitem{talagrandbk}
%
\begin{binproceedings}[vtex]
\bauthor{\bsnm{Talagrand},~\bfnm{Michel}\binits{M.}}
(\byear{1994}).
\btitle{Some remarks on the {Berg--Kesten} inequality}.
In \bbooktitle{Probability in {Banach} Spaces, 9 (Sandjberg, 1993)}.
\bseries{Progress in Probability}
\bvolume{35}
\bpages{293--297}.
\bpublisher{Birkh\"auser}, \baddress{Boston, MA}.
\bid{mr={1308525}}
\end{binproceedings}
%
\endbibitem

%b15 ###
\bibitem{vdBKBK}
%
\begin{barticle}[author]
\bauthor{\bparticle{van~den }\bsnm{Berg},~\bfnm{Jacob}\binits{J.}}
\AND
\bauthor{\bsnm{Kesten},~\bfnm{Harry}\binits{H.}}
(\byear{1985}).
\btitle{Inequalities with applications to percolation and reliability}.
\bjournal{J. Appl. Probab.}
\bvolume{22}
\bpages{556--569}.
\bid{mr={799280}}
\end{barticle}
%
\endbibitem

%b16 ###
\bibitem{wernerparkcity}
%
\begin{bincollection}[vtex]
\bauthor{\bsnm{Werner},~\bfnm{Wendelin}\binits{W.}}
(\byear{2009}).
\btitle{Lectures on two-dimensional critical percolation}.
In \bbooktitle{Statistical Mechanics}.
\bseries{IAS/Park City Math. Ser.}
\bvolume{16}
\bpages{297--360}.
\bpublisher{Amer. Math. Soc.}, \baddress{Providence, RI}.
\bid{mr={2523462}}
\end{bincollection}
%
\endbibitem

\end{thebibliography}
\end{document}